# Primality test and primes enumeration using odd numbers indexation

**WOLF Marc,** https://orcid.org/0000-0002-6518-9882

**WOLF François,** https://orcid.org/0000-0002-3330-6087

*Independent researchers;*

marc.wolf3@wanadoo.fr ; francois.wolf@dbmail.com

April 5th, 2020

## ABSTRACT

Odd numbers can be indexed by the map $k(n) = (n-3)/2$, $n \in 2\mathbb{N} + 3$. We first propose a basic primality test using this index function that was first introduced in [8]. Input size of operations is reduced which improves computational time by a constant. We then apply similar techniques to Atkin's prime-numbers sieve which uses modulus operations and finally to Pritchard's wheel sieve, in both case yielding similar results.

**Keywords**: odd number index, primality test, primes enumeration, Atkin sieve, composite odd numbers, wheel sieve.

### *CONTENTS*

# 1 Introduction

## 1.1 Primality test and prime enumeration

An odd number $n$ is prime when it is not divisible by any prime $p$ lower than or equal to $\sqrt{n}$. This basic primality test requires too much computational time for large integers. Faster and more efficient deterministic and probabilistic primality tests have been designed for large numbers [1]. A deterministic polynomial primality test was proposed by M. Agrawal, N. Kayal and N. Saxena in 2002 [2].

Enumeration of primes up to a given limit can be done by using a primality test but prime number sieves are preferred from a performance point of view. A sieve is a type of fast algorithm to find all primes up to a given number. There exists many such algorithms, from the simple Erastosthenes' sieve (invented more than 2000 years ago), to the wheel sieves of Paul Pritchard ([3], [4], [5]) and the sieve of Atkin [6]. In [7], Gabriel Paillard, Felipe Franca and Christian Lavault present another version of the wheel sieve and give an overview of all the existing prime-numbers sieves.

In theory, indices are a way to represent odd numbers. By adapting results from [8], we show how odd number indices may be used in applied mathematics. In the last part, we apply [8] to Pritchard's wheel sieve, which leads to a *dynamical* wheel sieve. Using the linear diophantine equation resolution method first introduced in [9], we introduce an original way of "turning the wheel".

## 1.2 Notation

We will use the following notations:

1. $I$ designates the set of odd integers greater than 1, i.e.:

$$I = \{N_k = 2k + 3 | k \in \mathbb{N}\};$$

2. $P$ the set of prime numbers, $P_n$ the set of prime numbers not greater than $n$;

3. $C$ the set of composite odd integers, i.e.:

$$C = I \backslash P = \{N_k \in I | \exists (a, b) \in I, N_k = ab\}$$

The function $f: k \in \mathbb{N} \mapsto N_k \in I$ is bijective. The inverse function is $f^{-1}: N_k \in I \mapsto k = \frac{N_k - 3}{2}$. $k = f^{-1}(N_k)$ is the index of $N_k$. The preimage of $C$ is denoted by $W$:

$$W = f^{-1}(C) = \{k \in \mathbb{N} \mid N_k \in C\}$$

4. For $x$ and $y$ two integers, we denote by $x \bmod y$ the remainder of the Euclidean division of $x$ by $y$, which belongs to $[\![0, y-1]\!]$.

5. $N_1$ and $N_2$ are the subsets of $I$ given by:

$$N_1 = \{N_k \in I \mid N_k \bmod 4 = 1\}$$
$$N_2 = \{N_k \in I \mid N_k \bmod 4 = 3\}$$

Similarly:

$$C_1 = N_1 \cap C$$
$$C_2 = N_2 \cap C$$

Finally, $S_1$ and $S_2$ designate the set of indices corresponding to elements of $C_1$ and $C_2$ respectively, i.e. $S_1 = f^{-1}(C_1)$ and $S_2 = f^{-1}(C_2)$.

# 2 Basic primality test and primes enumeration

## 2.1 Two families of infinite sequences with arithmetic difference

[8] shows that $W$ is the union of two families of finite sequences with arithmetic difference. Actually proposition 2-5 says that any composite odd number $N_k \in C$ can be written as a difference of two squares, and more precisely that there exists $j \in \mathbb{N}$ and $x \in [\![0, j]\!]$ such that:

$$\begin{cases} \textbf{(1) } N_k \in C_1 \Rightarrow N_k = (2j+3)^2 - (2x)^2, \\ \textbf{(2) } N_k \in C_2 \Rightarrow N_k = (2j+4)^2 - (2x+1)^2 \end{cases}$$

**Corollary 2-1**: Let $k_j(n) = (2j+3)n + j$. One has:

$$W = S_1 \cup S_2$$

and:

$$\begin{aligned} S_1 &= \{k_i(x) = k_i(i+1) + 2(2i+3)x; \ i \in \mathbb{N}, x \in \mathbb{N}\} \\ S_2 &= \{k_i(x) = k_i(i+2) + 2(2i+3)x; i \in \mathbb{N}, x \in \mathbb{N}\} \end{aligned}$$

Thus $W$ is the union of two families of infinite arithmetic sequences. The indices $k_i(i+1)$ of first type reference points (or remarkable points, see[8]) are the initial terms of sequences ranging in $S_1$. Similarly, the indices $k_i(i+2)$ of second type reference points are the initial terms of sequences ranging in $S_2$.

*Proof*: We substitute $j$ by $i+x$ in relations **(1)** and **(2)**:

$$(2j+3)^2 - (2x)^2 = (2i+2x+3)^2 - (2x)^2 = (2i+3)(2i+4x+3)$$
$$= 2[k_i(i+1) + 2(2i+3)x] + 3$$

and similarly:

$$(2j+4)^2 - (2x+1)^2 = (2i+2x+4)^2 - (2x+1)^2 = (2i+3)(2i+4x+5)$$
$$= 2(2i+3)(i+2x+2) + 2i + 3 = 2[k_i(i+2) + 2(2i+3)x] + 3$$

**Proposition 2-1**: For any $N_k \in C$ there exists $X \in P, X \leq \sqrt{N_k}$ and $x \in \mathbb{N}$ such that:

$$N_k \in C_1 \Rightarrow N_k = X(X + 4x)$$

$$N_k \in C_2 \Rightarrow N_k = X(X + 4x + 2)$$

Thus, writing $X = 2i + 3$, we get:

$$W = S_1' \cup S_2'$$

where:

$$S_1' = \{k_i(x) = k_i(i+1) + 2(2i+3)x; \; i \in \mathbb{N} \setminus W, x \in \mathbb{N}\}$$
$$S_2' = \{k_i(x) = k_i(i+2) + 2(2i+3)x; \; i \in \mathbb{N} \setminus W, x \in \mathbb{N}\}$$

*Proof*: Take $X$ the smallest prime dividing $N_k \in C$. Thus $X \in P_{\sqrt{N_k}}$ and if $Y = \frac{N_k}{X}$ then $Y \geq X$ and $Y - X$ is even, and we can write it either $4x$ or $4x + 2$. These two cases clearly correspond respectively to $N_k \in C_1$ and $N_k \in C_2$. Thus the index $k$ can be decomposed as in corollary 2-1, but with $i$ the index of a prime number, hence in $\mathbb{N} \setminus W$.

## 2.2 Basic primality test

In this section, we describe a basic primality test using the previous infinite sequences.

**Definition 2-2**: For any $p = 2i + 3 \in P$ and $N \in I$ we let:

1- $A(N, p) = N - p^2$ and $f_A(p) = p^2$.

2- $B(N, p) = N - p(p + 2)$ and $f_B(p) = p(p + 2)$.

**Proposition 2-2**: $N \in N_1$ is a prime number when:

$$\forall p = 2i + 3 \in P_{\sqrt{N}}, \frac{A(N,p)}{4} \bmod p \neq 0$$

$N \in N_2$ is a prime number when:

$$\forall p = 2i + 3 \in P_{\sqrt{N}}, \frac{B(N,p)}{4} \bmod p \neq 0$$

*Proof*: This follows from the fact that $A(N, p) \bmod p = N \bmod p$ and likewise for $B(N, p)$.

**Remark 2-2:** In order to reduce computation of $A(N, p)$ and $B(N, p)$ for two consecutive prime numbers, we only decrement the value.

More precisely, if $p < p'$ are two primes, we let $\alpha(p, p') = p' - p$ and we compute:

$$\begin{cases} \Delta A(N, p, p') = A(N, p) - A(N, p') = \alpha(\alpha + 2p) \\ \Delta B(N, p, p') = B(N, p) - B(N, p') = \Delta A(N, p, p') + 2\alpha \end{cases}$$

These two expressions are independent of $N$.

### 2.3 Primality test with indices

We adapt here the results of the previous section with indices.

**<u>Definition 2-3</u>**: For any $i$ index of a prime number $p \in P$ and $k \in \mathbb{N}$, we let:

1- $A'(k, i) = (k - 3)/2 - i(i + 3)$, $f'_A(i) = i(i + 3)$, $g'_A(k) = (k - 3)/2$

2- $B'(k, i) = (k - 6)/2 - i(i + 4)$ and $f'_B(i) = i(i + 4)$, $g'_B(k) = (k - 6)/2$

**<u>Proposition 2-3</u>**: $k \in S_1$ is a prime number index when:

$$\forall p = 2i + 3 \in P_{\sqrt{2k+3}}, A'(k, i) \bmod p \neq 0$$

$k \in S_2$ is a prime number index when:

$$\forall p = 2i + 3 \in P_{\sqrt{2k+3}}, B'(k, i) \bmod p \neq 0$$

<u>*Proof*</u>: This follows from proposition 2-2 and definition 2-2 because if we let $N = 2k + 3$ then $A'(k, i) = \frac{A(N,p)}{4}$ and $B'(k, i) = \frac{B(N,p)}{4}$.

**<u>Remark 2-3:</u>** In order to reduce computation of $A'(k, i)$ and $B'(k, i)$ for two consecutive prime number indices, we only decrement their values.

More precisely, if $i < i'$ are two prime indices we let $\alpha'(i, i') = i' - i$ and we compute:

$$\Delta A'(k, i, i') = A'(k, i) - A'(k, i') = \alpha'(\alpha' + 2i + 3)$$

$$\Delta B'(k, i, i') = B'(k, i) - B'(k, i') = \Delta A'(k, i, i') + \alpha'$$

These two expressions are independent of $k$.

### 2.4 First algorithms of prime enumeration

In this section, we present prime enumeration algorithms based on propostion 2-2 and 2-3. The first one manipulates numbers and the second one indices.

#### 2.4.1 Primality test using numbers

This first algorithm named **_PrimeEnumeration_** consists in two functions:

- The main function which determines primes in up to $N_{Max}$ and returns them in a list, along with its size.
- An auxiliary function which returns whether a number $N$ is prime, based on precomputed list of primes and values of $\Delta A$ and $\Delta B$. It is called **_LocalTest._** It is also in charge of updating the lists $\Delta A$ and $\Delta B$ if needed.

Three zero-based lists are used and built recursively in this algorithm: the list of primes itself $L_p$, and the lists of values for $\Delta A$ and $\Delta B$ respective to $L_p$ (remember it is independent from $N$). Only numbers which are not multiples of 2 and 3 are tested. Thus we restrict to $N = 6m + 1$ and $N = 6m + 5$. The congruence of $N$ modulo 4 depends on the parity of $m$, i.e. when $m$ is even, $N \bmod 4 = 1$ and when $m$ is odd, $N \bmod 4 = 3$.

---

**Algorithm 2-4-1a Function *PrimeEnumeration($N_{Max}$)*:** $N_{Max}$ is an odd integer such that $N_{Max} \geq 7$. This function returns the list of primes up to $N_{Max}$ and its size.

---

***First step : intialisation of variables***

$L_p \leftarrow \{5\}$        → List of primes from 5, initialized with one element

$i_l \leftarrow 1$        → Size of the list $L_p$

       **→ About the next two lists, see the remark 2-2**

$\Delta A \leftarrow \{16\}$        → $\Delta A(N, 3,5) = 2 \times (2 + 6) = 16$

$\Delta B \leftarrow \{20\}$        → $\Delta B(N, 3,5) = \Delta A(N, 3,5) + 2 \times 2 = 20$

$i_{r1} \rightarrow 0$

$Cap1 \leftarrow 25$

$i_{r2} \rightarrow 0$

$Cap2 \leftarrow 35$

***Second step : iteration***

$(m, N) \leftarrow (1,7)$

$ModEqOne \leftarrow$ **False**        → $m = 1$ so $(6m + 1) \bmod 4 = 3$

**While** $N \leq N_{Max}$ **Do**        → Loop to get odd primes in range $[\![5, N_{Max}]\!]$

  **If** *LocalTest*$(N, L_p, \Delta A, \Delta B, i_{r1}, Cap1, i_{r2}, Cap2, ModEqOne)$ **Do**

    $L_p(i_l) \leftarrow N$

    $i_l \leftarrow i_l + 1$

  **End If**

  $N \leftarrow 6m + 5$

  **If** $N \leq N_{Max}$ **And** *LocalTest*$(N, L_p, \Delta A, \Delta B, i_{r1}, Cap1, i_{r2}, Cap2, ModEqOne)$ **Do**

    $L_p(i_l) \leftarrow N$

    $i_l \leftarrow i_l + 1$

  **End If**

  $m \leftarrow m + 1$

$N \leftarrow 6m + 1$

$ModEqOne \leftarrow !\,ModEqOne$  → Switch the boolean value

**End While**

**Return** $(\{2,3\} + L_p, i_l + 2)$  → Return the list of primes and the number of primes.

---

**Algorithm 2-4-1b Function** *LocalTest* $(N, L_p, \Delta A, \Delta B, i_{r1}, Cap1, i_{r2}, Cap2, ModEqOne)$**:** $N$ is an odd integer. $i_{root}$ stands for $i_{r1}$ or $i_{r2}$ depending on $ModEqOne$. This function decides whether for all $p \in L_p[0 \dots i_{root}]$, $A(N,p)/4$ or $B(N,p)/4$ is not divisible by $p$. It will also potentially update $\Delta A$, $\Delta B$, $i_{r1}$, $i_{r2}$, $Cap1$ and $Cap2$ which must be passed by reference.

---

***First step : intialisation of variables***

$A \leftarrow 9$  → stands for $f_A(3) = 3^2$

$B \leftarrow 15$  → stands for $f_B(3) = 3 \times 5$

**If** $ModEqOne$ **Do**  → initiate references that might be updated

$i_{root} \leftarrow i_{r1}$

$Cap \leftarrow Cap1$

$\Delta \leftarrow \Delta A$

**Else**

$i_{root} \leftarrow i_{r2}$

$Cap \leftarrow Cap2$

$\Delta = \Delta B$

**End If**

**If** $N = Cap$ **Do**

  **Return False**  → The cap is a composite number

**End If**

**If** $N > Cap$ **Do**  → update references because we always want $N \leq Cap$

$i_{root} \leftarrow i_{root} + 1$

$\alpha \leftarrow \left( L_p(i_{root}) - L_p(i_{root} - 1) \right)$

  **If** $ModEqOne$ **Do**

  $\Delta(i_{root}) \leftarrow \alpha(\alpha + 2L_p(i_{root} - 1))$  → $\Delta A$

  **Else**

  $\Delta(i_{root}) \leftarrow \Delta A(i_{root}) + 2\alpha$  → $\Delta B$, using $\Delta A$ which must already be updated

**End If**

$$Cap \leftarrow Cap + \Delta(i_{root})$$

**End If**

***Second step : iteration***

**If** $ModEqOne$ **Do**

$$N \leftarrow N - A$$

**Else**

$$N \leftarrow N - B$$

**End If**

$i \leftarrow 0$

**While** $i \leq i_{root}$ **Do**          → Iteration at most up to $i = i_{root}$

$\quad N \leftarrow N - \Delta(i)$

$\quad$ **If** $(N/4) \bmod L_p(i) = 0$ **Do**    → $N$ is a multiple of 4, division by 4 can be done bitwise

$\quad\quad$ **Return False**           → Test is negative

$\quad$ **End If**

$\quad i \leftarrow i + 1$

**End While**

**Return True**               → Test is positive

---

### 2.4.2 Primality test using infinite sequences and indices

This second algorithm ***IndexPrimeEnumeration*** also consists in two functions, mirroring the previous algorithm:

> ➢ The main function which determines primes up to $N_{Max}$ and returns them in a list along with its size.
> ➢ An auxiliary function which returns whether a number $N$ is prime based on precomputed list of primes and values of $\Delta A'$ and $\Delta B'$. It is called ***LocalTest.***

Four zero-based lists are used and built recursively: the list of primes $L_p$, the corresponding indices $IL_p$ (indices of primes), and the lists $\Delta A'$ and $\Delta B'$ respective to $L_p$.

Only numbers which are not multiple of 2 and 3 are tested, i.e. indices of the form $k = 3m - 1$ and $k = 3m + 1$.

**Remark 2-4-2**: To avoid any division in the computation of $A'$ and $B'$ we will write $m = 2t + 1$ or $2t + 2$.

---

**Algorithm 2-4-2a Function *IndexPrimeEnumeration($N_{Max}$)***: $N_{Max}$ is an odd integer such that $N_{Max} \geq 7$. This function returns the list of primes up to $N_{Max}$ and its size.

---

***First step : intialisation of variables***

$L_p \leftarrow \{5\}$              → List of primes from 5, initialized with one element

$IL_p \leftarrow \{1\}$              → List of index of primes

$i_l \leftarrow 1$              → Size of the two lists $L_p$ and $IL_p$

→ **About the next two lists, see the remark 2-3**

$\Delta A' \leftarrow \{4\}$              → $\Delta A'(k, 0,1) = 1 \times (1 + 3) = 4$

$\Delta B' \leftarrow \{5\}$              → $\Delta B'(k, 0,1) = \Delta A'(k, 0,1) + 1 = 5$

$k_{Max} \leftarrow (N_{Max} - 3)/2$

$i_{r1} \rightarrow 0$

$Cap1 \leftarrow 11$

$i_{r2} \rightarrow 0$

$Cap2 \leftarrow 16$

***Second step : iteration***

$(t, k, g') \leftarrow (0, 2, -2)$              → $k$ starts at $3(2t + 1) - 1$, $g'$ stands for $g'_A(k)$ or $g'_B(k)$

**While** $k \leq k_{Max}$ **Do**              → Loop to get odd prime indices in range $[\![1, k_{Max}]\!]$

  **If** *LocalTest($g', k, L_p, IL_p, \Delta A', \Delta B', i_{r2}, Cap2$,**False**)* **Do**

  $IL_p(i_l) \leftarrow k$

  $L_p(i_l) \leftarrow 2k + 3$

  $i_l \leftarrow i_l + 1$

  **End If**

  $k \leftarrow k + 2$              → $k = 3(2t + 1) + 1$

  $g' \leftarrow g' + 1$

  **If** $k \leq k_{Max}$ **And** *LocalTest($g', k, L_p, IL_p, \Delta A', \Delta B', i_{r2}, Cap2$,**False**)* **Do**

  $IL_p(i_l) \leftarrow k$

  $L_p(i_l) \leftarrow 2k + 3$

$i_l \leftarrow i_l + 1$

**End If**

$k \leftarrow k + 1$                       → $k = 3(2t + 2) - 1$

$g' \leftarrow g' + 2$

**If** $k \leq k_{Max}$ **And** $LocalTest(g', k, L_p, IL_p, \Delta A', \Delta B', i_{r1}, Cap1,$ **True**) **Do**

   $IL_p(i_l) \leftarrow k$

   $L_p(i_l) \leftarrow 2k + 3$

   $i_l \leftarrow i_l + 1$

**End If**

$k \leftarrow k + 2$                       → $k = 3(2t + 2) + 1$

$g' \leftarrow g' + 1$

**If** $k \leq k_{Max}$ **And** $LocalTest(g', k, L_p, IL_p, \Delta A', \Delta B', i_{r1}, Cap1,$ **True**) **Do**

   $IL_p(i_l) \leftarrow k$

   $L_p(i_l) \leftarrow 2k + 3$

   $i_l \leftarrow i_l + 1$

**End If**

$t \leftarrow t + 1$                       → We do not use $t$ but keep it for the sake of readability

$k \leftarrow k + 1$                       → $k = 3(2t + 1) - 1$

$g' \leftarrow g' - 1$

**End While**

**Return** $(\{2,3\} + L_p, i_l + 2)$       → Return the list of primes and the number of primes.

---

**Algorithm 2-4-2b Function** $LocalTest(g', k, L_p, IL_p, \Delta A', \Delta B', i_{root}, Cap, ModEqOne)$**:** $g'$ stands for $g'_A(k)$ or $g'_B(k)$ depending on $ModEqOne$. This function decides whether for all $p \in L_p[0 \dots i_{root}]$, $A'(k, i)$ or $B'(k, i)$ is coprime with $p$.

---

***First step : intialisation of variables***

**If** $ModEqOne$ **Do**                 → initiate references that might be updated

   $\Delta \leftarrow \Delta A'$

**Else**

   $\Delta = \Delta B'$

**End If**

**If** $k = Cap$ **Do**

  **Return False**           → The cap is the index of a composite number

**End If**

**If** $k > Cap$ **Do**          → update references because we always want $k \leq Cap$

$i_{root} \leftarrow i_{root} + 1$

$\alpha \leftarrow \left( IL_p(i_{root}) - IL_p(i_{root} - 1) \right)$

**If** $ModEqOne$ **Do**

  $\Delta(i_{root}) \leftarrow \alpha(\alpha + L_p(i_{root} - 1))$     → $\Delta A'$

**Else**

  $\Delta(i_{root}) \leftarrow \Delta A'(i_{root}) + \alpha$     → $\Delta B'$, using $\Delta A'$ which must already be updated

**End If**

$Cap \leftarrow Cap + \Delta(i_{root})$

**End If**

*Second step : iteration*

$R \leftarrow g'$

$i \leftarrow 0$

**While** $i \leq i_{root}$ **Do**       → Iteration at most up to $i = i_{root}$

$R \leftarrow R - \Delta(i)$

**If** $R \bmod L_p(i) = 0$ **Do**

  **Return False**        → Test is negative

  **End If**

$i \leftarrow i + 1$

**End While**

**Return True**        → Test is positive

---

### 2.5 Perfomance of the algorithms

In this section, we present the performance of the previous two algorithms of prime enumeration. We first give a theoretical complexity, followed by empirical results.

<u>**Proposition 2-5:**</u> Time complexity (in terms of number of arithmetic operations) and space complexity are the same for both *PrimeEnumeration* and *IndexPrimeEnumeration* algorithms*.*

Time complexity is:

$$O\left(\frac{(N_{Max})^{\frac{3}{2}}}{\ln(N_{Max})}\right).$$

Space complexity is:

$$O\left(\frac{\sqrt{N_{Max}}}{\ln(N_{Max})}\right)$$

<u>*Proof:*</u> Any number $n$'s primality is tested with primes in $[\![5, \sqrt{n}]\!]$, in $O(1)$ operations. There are $\pi(\sqrt{n}) - 2 \sim \frac{\sqrt{n}}{\ln(\sqrt{n})} = O\left(\frac{\sqrt{n}}{\ln(n)}\right)$ such primes. We loop over range $[\![7, N_{Max}]\!]$, time complexity is thus $\sum_{t=7}^{N_{Max}} O\left(\frac{\sqrt{t}}{\ln(t)}\right) = O\left(\frac{(N_{Max})^{\frac{3}{2}}}{\ln(N_{Max})}\right)$ (actually we skip two thirds of the terms in this sum by not testing multiples of 2 and 3, but complexity remains $O\left(\frac{(N_{Max})^{\frac{3}{2}}}{\ln(N_{Max})}\right)$ albeit with smaller constant.

The space complexity is related to the lists we keep in memory, which are at most of size $\pi(N_{Max})$. This space complexity is $O\left(\frac{\sqrt{N_{Max}}}{\ln(N_{Max})}\right)$.

Both algorithms have been implemented in Visual Studio C++ 2012. We measured execution time for various values of $N_{Max}$ and produced a regression using Maple 2017.3. Details of the Maple options used to get the regression are given in appendix 8.1.

On the graph 2-5 below, we represent the computation time in seconds for both algorithms. Curve $T_1$ corresponds to the algorithm ***PrimeEnumeration*** and curve $T_2$ to ***IndexPrimeEnumeration***. The correlation coefficient R of each curve is given on the graph. We observe that computation time of both algorithms is consistent with theoretical complexity, although exponent is a bit smaller than 1.5.

**<u>Graph 2-5</u>: computation time $T\ (N_{Max})$ in seconds for both algorithms (Prime enumeration)**

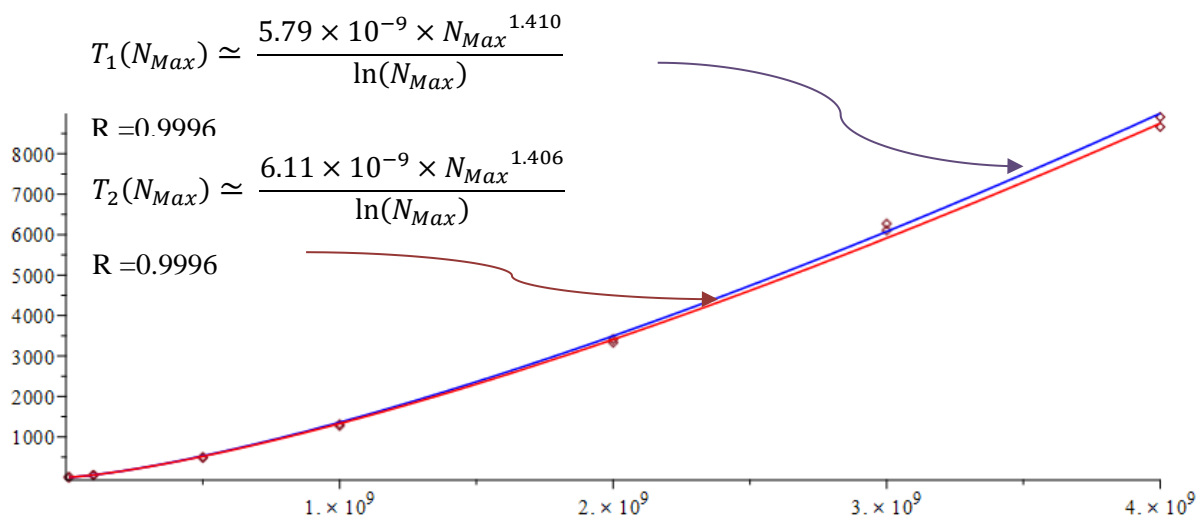

$$T_1(N_{Max}) \simeq \frac{5.79 \times 10^{-9} \times N_{Max}^{\,1.410}}{\ln(N_{Max})}$$

R =0.9996

$$T_2(N_{Max}) \simeq \frac{6.11 \times 10^{-9} \times N_{Max}^{\,1.406}}{\ln(N_{Max})}$$

R =0.9996

Both algorithms **PrimeEnumeration** and **IndexPrimeEnumeration** have the same number of modulo operations. But the computation of the input of modulus operations is done with larger inputs for the former than for the latter, which allows to marginally save time for large values of $N_{Max}$.

# 3 The sieve of Atkin

The sieve of Atkin [6] is a modern and efficient algorithm for primes enumeration. We present two algorithms based on it, one using numbers and the other indices. Both are based on the version which has a complexity $O(N_{Max})$ in time and space. Modified versions achieve up to $O\left(\frac{N_{Max}}{\ln\ln(N_{Max})}\right)$ in time and $O\left(N_{Max}^{\frac{1}{2}+o(1)}\right)$ in space.

## 3.1 Atkin algorithm

This algorithm is based on the following three results from [6].

**Proposition 3-1** Let $n > 3$ be a square-free integer. Then $n$ is prime if and only if one of the three following conditions is true:

    a. $n \in 1 + 4\mathbb{N}$ and there is an odd number of solutions to $n = 4x^2 + y^2, (x, y) \in \mathbb{N}^2$,
    b. $n \in 7 + 12\mathbb{N}$ and there is an odd number of solutions to $n = 3x^2 + y^2, (x, y) \in \mathbb{N}^2$,
    c. $n \in 11 + 12\mathbb{N}$ and there is an odd number of solutions to $n = 3x^2 - y^2, x > y, (x, y) \in \mathbb{N}^2$.

We observe that the first congruence condition on $n$ can also be replaced by $n \in 1 + 12\mathbb{N}$ or $n \in 5 + 12\mathbb{N}$. We also observe the following for an odd integer $n$:

- If $n = 4x^2 + y^2$, $y$ must be odd.
- If $n = 3x^2 + y^2$ or $n = 3x^2 - y^2$, $x$ and $y$ must have opposite parity.

Furthermore if $n$ is square-free, $x$ and $y$ must be in $\mathbb{N}^*$, with $x < \sqrt{n/2}$ and $y < \sqrt{n}$.

**Remark 3-1** We can compute the remainder modulo 12 of $ax^2 + by^2$ depending on remainders modulo 12 of $x$ and $y$. This gives us the different cases to check in Atkin sieve. We present them in table 3-1, noting that there is no case for $y \bmod 12 = 0$ and $y \bmod 12 = 6$.

**Table 3-1: Atkin sieve cases depending on remainders modulo 12 of $x$ and $y$.**

| $x\backslash y$ | 1 | 2 | 3 | 4 | 5 | 7 | 8 | 9 | 10 | 11 |
|---|---|---|---|---|---|---|---|---|---|---|
| 0 | $4x^2+y^2$ $3x^2-y^2$ | | | | $4x^2+y^2$ $3x^2-y^2$ | $4x^2+y^2$ $3x^2-y^2$ | | | | $4x^2+y^2$ $3x^2-y^2$ |
| 1 | $4x^2+y^2$ | $3x^2+y^2$ $3x^2-y^2$ | $4x^2+y^2$ | $3x^2+y^2$ $3x^2-y^2$ | $4x^2+y^2$ | $4x^2+y^2$ | $3x^2+y^2$ $3x^2-y^2$ | $4x^2+y^2$ | $3x^2+y^2$ $3x^2-y^2$ | $4x^2+y^2$ |
| 2 | $4x^2+y^2$ $3x^2-y^2$ | | $4x^2+y^2$ | | $4x^2+y^2$ $3x^2-y^2$ | $4x^2+y^2$ $3x^2-y^2$ | | $4x^2+y^2$ | | $4x^2+y^2$ $3x^2-y^2$ |
| 3 | $4x^2+y^2$ | $3x^2+y^2$ $3x^2-y^2$ | | $3x^2+y^2$ $3x^2-y^2$ | $4x^2+y^2$ | $4x^2+y^2$ | $3x^2+y^2$ $3x^2-y^2$ | | $3x^2+y^2$ $3x^2-y^2$ | $4x^2+y^2$ |
| 4 | $4x^2+y^2$ $3x^2-y^2$ | | $4x^2+y^2$ | | $4x^2+y^2$ $3x^2-y^2$ | $4x^2+y^2$ $3x^2-y^2$ | | $4x^2+y^2$ | | $4x^2+y^2$ $3x^2-y^2$ |

| 5 | $4x^2+y^2$ | $3x^2+y^2$ $3x^2-y^2$ | $4x^2+y^2$ | $3x^2+y^2$ $3x^2-y^2$ | $4x^2+y^2$ | $4x^2+y^2$ | $3x^2+y^2$ $3x^2-y^2$ | $4x^2+y^2$ | $3x^2+y^2$ $3x^2-y^2$ | $4x^2+y^2$ |
|---|---|---|---|---|---|---|---|---|---|---|
| 6 | $4x^2+y^2$ $3x^2-y^2$ | | | | $4x^2+y^2$ $3x^2-y^2$ | $4x^2+y^2$ $3x^2-y^2$ | | | | $4x^2+y^2$ $3x^2-y^2$ |
| 7 | $4x^2+y^2$ | $3x^2+y^2$ $3x^2-y^2$ | $4x^2+y^2$ | $3x^2+y^2$ $3x^2-y^2$ | $4x^2+y^2$ | $4x^2+y^2$ | $3x^2+y^2$ $3x^2-y^2$ | $4x^2+y^2$ | $3x^2+y^2$ $3x^2-y^2$ | $4x^2+y^2$ |
| 8 | $4x^2+y^2$ $3x^2-y^2$ | | $4x^2+y^2$ | | $4x^2+y^2$ $3x^2-y^2$ | $4x^2+y^2$ $3x^2-y^2$ | | $4x^2+y^2$ | | $4x^2+y^2$ $3x^2-y^2$ |
| 9 | $4x^2+y^2$ | $3x^2+y^2$ $3x^2-y^2$ | | $3x^2+y^2$ $3x^2-y^2$ | $4x^2+y^2$ | $4x^2+y^2$ | $3x^2+y^2$ $3x^2-y^2$ | | $3x^2+y^2$ $3x^2-y^2$ | $4x^2+y^2$ |
| 10 | $4x^2+y^2$ $3x^2-y^2$ | | $4x^2+y^2$ | | $4x^2+y^2$ $3x^2-y^2$ | $4x^2+y^2$ $3x^2-y^2$ | | $4x^2+y^2$ | | $4x^2+y^2$ $3x^2-y^2$ |
| 11 | $4x^2+y^2$ | $3x^2+y^2$ $3x^2-y^2$ | $4x^2+y^2$ | $3x^2+y^2$ $3x^2-y^2$ | $4x^2+y^2$ | $4x^2+y^2$ | $3x^2+y^2$ $3x^2-y^2$ | $4x^2+y^2$ | $3x^2+y^2$ $3x^2-y^2$ | $4x^2+y^2$ |

We could run the sieve looping through 12x12 blocks of $(x, y)$ according to this table, but for readability we do not implement this optimization in the algorithms below. We note however that this would save all the modulo operations.

---

**Algorithm 3-1 *SieveOfAtkin($N_{Max}$)*:** $N_{Max} > 3$ is an integer. This function returns the list of all prime numbers less than $N_{Max}$.

---

***First step : intialisation of variables***

$L_p \leftarrow \{2,3\}$            $\rightarrow$ Dynamic list of odd primes

$i_l \leftarrow 2$              $\rightarrow$ Number of primes in the list

Sieve[$N_{Max}$] $\leftarrow$ {**False**, …,**False**} $\rightarrow$ Array of $N_{Max}$ entries all initialized to **False**

$x_{max} \leftarrow \lceil \sqrt{N_{Max}/2} \rceil - 1$     $\rightarrow$ Bound for $x$

$y_{max} \leftarrow \lceil \sqrt{N_{Max}} \rceil - 1$       $\rightarrow$ Bound for $y$

***Second step : iteration for first case***

**For** $x = 1$ **To** $x_{max}$

  **For** $y = 1$ **To** $y_{max}$ **Step** $2$     $\rightarrow$ $y$ must be odd

  $n \leftarrow 4x^2 + y^2$

  **If** $n < N_{Max}$ **And** ($n \bmod 12 = 1$ **Or** $n \bmod 12 = 5$) **Do**

  Sieve[$n$] $\leftarrow$ !Sieve[$n$]     $\rightarrow$ Switch the boolean value Sieve[$n$]

**End If**

    **End For**

**End For**

*Third step : iteration for second and third cases*

**For** $x = 1$ **To** $x_{max}$ **Step** $2$

  **For** $y = 2$ **To** $y_{max}$ **Step** $2$       → case where $x$ is odd and $y$ even

   $n \leftarrow 3x^2 + y^2$

   **If** $n < N_{Max}$ **And** $(n \bmod 12 = 7)$ **Do**

    Sieve$[n] \leftarrow$ !Sieve$[n]$

   **End If**

   **If** $x > y$ **Do**

    $n \leftarrow 3x^2 - y^2$

    **If** $n < N_{Max}$ **And** $(n \bmod 12 = 11)$ **Do**

     Sieve$[n] \leftarrow$ !Sieve$[n]$

    **End If**

   **End If**

  **End For**

**End For**

**For** $x = 2$ **To** $x_{max}$ **Step** $2$

  **For** $y = 1$ **To** $y_{max}$ **Step** $2$       → case where $x$ is even and $y$ is odd

   $n \leftarrow 3x^2 + y^2$

   **If** $n < N_{Max}$ **And** $(n \bmod 12 = 7)$ **Do**

    Sieve$[n] \leftarrow$ !Sieve$[n]$

   **End If**

   **If** $x > y$ **Do**

    $n \leftarrow 3x^2 - y^2$

    **If** $n < N_{Max}$ **And** $(n \bmod 12 = 11)$ **Do**

     Sieve$[n] \leftarrow$ !Sieve$[n]$

    **End If**

   **End If**

  **End For**

**End For**

***Fourth step : remove multiples of prime squares***

**For** $n = 5$ **To** $y_{max}$ **Step** $2$      → multiples of 2 and 3 are ignored by the previous iterations

  **If** Sieve[$n$] **Do**

   **For** $i = n^2$ **To** $N_{Max} - 1$ **Step** $2n^2$

   Sieve[$i$] ← **False**

   **End For**

  **End If**

**End For**

***Last step : return list of primes from the sieve***

**For** $n = 5$ **To** $N_{Max} - 1$ **Step** $2$

  **If** Sieve[$n$] **Do**

   $L_p(i_l) \leftarrow n$

   $i_l \leftarrow i_l + 1$

  **End If**

**End For**

**Return** $(L_p, i_l)$

---

### 3.2 Atkin algorithm with indices

We can rewrite proposition 3-1 as:

**Corollary 3-2:** $k$ is the index of a prime number if and only if $2k + 3$ is square-free and one of the three following conditions is true:

   a. $k \in (1 + 6\mathbb{N}) \cup (5 + 6\mathbb{N})$ and there is an odd number of solutions to $k = 2x^2 + \frac{y^2 - 3}{2}$,

   b. $k \in 2 + 6\mathbb{N}$ and there is an odd number of solutions to $k = \frac{3x^2 + y^2 - 3}{2}$,

   c. $k \in 4 + 6\mathbb{N}$ and there is an odd number of solutions to $k = \frac{3x^2 - y^2 - 3}{2}$ with $y < x$.

The relationships presented in the following remark are used in the next algorithm.

**Remark 3-2:** For the fourth step (square multiples elimination), we note that if $n = 2k + 3$, the index of $n^2$ is $2k^2 + 6k + 3$ and that the step of $2n^2$ translates into a step of $n^2 = (2k + 3)^2$ for indices.

---

**Algorithm 3-2** *IndexSieveOfAtkin($N_{Max}$)*: $N_{Max} > 3$ is an odd integer. This function returns the list of all prime numbers less than $N_{Max}$.

---

***First step : intialisation of variables***

$L_p \leftarrow \{2,3\}$             → Dynamic list of primes

$i_l \leftarrow 2$               → Number of primes in the list

$k_{Max} \leftarrow (N_{Max} - 3)/2$       → Index of $N_{Max}$

$\text{Sieve}[k_{Max}] \leftarrow \{\textbf{False}, \dots ,\textbf{False}\}$       → Array of $k_{Max}$ entries all initialized to **False**

$\boldsymbol{x_{max}} \leftarrow \lceil \sqrt{N_{Max}/2} \rceil - 1$      → Bound for $x$

$\boldsymbol{y_{max}} \leftarrow \lceil \sqrt{N_{Max}} \rceil - 1$       → Bound for $y$

***Second step : iteration for first case***

**For** $x = 1$ **To** $x_{max}$

  **For** $y = 1$ **To** $y_{max}$ **Step** $2$        → $y$ must be odd

    $k \leftarrow 2x^2 + \dfrac{y^2 - 3}{2}$

    **If** $k < k_{Max}$ **And** ($k \bmod 6 = 1$ **Or** $k \bmod 6 = 5$) **Do**

      $\text{Sieve}[n] \leftarrow\ !\text{Sieve}[n]$        → Switch the boolean value Sieve[$n$]

    **End If**

  **End For**

**End For**

***Third step : iteration for second and third cases***

**For** $x = 1$ **To** $x_{max}$ **Step** $2$

  **For** $y = 2$ **To** $y_{max}$ **Step** $2$        → case where $x$ is odd and $y$ even

    $k \leftarrow \dfrac{3x^2 + y^2 - 3}{2}$

    **If** $k < k_{Max}$ **And** ($k \bmod 6 = 2$) **Do**

      $\text{Sieve}[n] \leftarrow\ !\text{Sieve}[n]$

    **End If**

    **If** $x > y$ **Do**

      $k \leftarrow \dfrac{3x^2 - y^2 - 3}{2}$

      **If** $k < N_{Max}$ **And** ($k \bmod 6 = 4$) **Do**

        $\text{Sieve}[n] \leftarrow\ !\text{Sieve}[n]$

      **End If**

    **End If**

  **End For**

**End For**

**For** $x = 2$ **To** $x_{max}$ **Step** $2$

  **For** $y = 1$ **To** $y_{max}$ **Step** $2$         $\rightarrow$ case where $x$ is even and $y$ is odd

    $k \leftarrow \frac{3x^2 + y^2 - 3}{2}$

    **If** $k < k_{Max}$ **And** $(k \bmod 6 = 2)$ **Do**

      Sieve$[n] \leftarrow$ !Sieve$[n]$

    **End If**

    **If** $x > y$ **Do**

      $k \leftarrow \frac{3x^2 - y^2 - 3}{2}$

      **If** $k < k_{Max}$ **And** $(k \bmod 6 = 4)$ **Do**

        Sieve$[n] \leftarrow$ !Sieve$[n]$

      **End If**

    **End If**

  **End For**

**End For**

*Fourth step : remove multiples of prime squares*

**For** $k = 1$ **To** $\frac{y_{max} - 3}{2}$         $\rightarrow$ multiples of 3 are ignored by the previous iterations

  **If** Sieve$[k]$ **Do**

    **For** $i = 2k^2 + 6k + 3$ **To** $k_{Max} - 1$ **Step** $(2k + 3)^2$

    Sieve$[i] \leftarrow$ **False**

    **End For**

  **End If**

**End For**

*Last step : return list of primes from the sieve*

**For** $k = 1$ **To** $k_{Max} - 1$

  **If** Sieve$[k]$ **Do**

    $L_p(i_l) \leftarrow 2k + 3$

    $i_l \leftarrow i_l + 1$

  **End If**

**End For**

**Return** $(L_p, i_l)$

### 3.3 Performance of algorithms

In this section, we discuss theoretical complexity and present our results with the two algorithms implementing the sieve of Atkin.

The reference algorithm ***SieveOfAtkin*** has less operations index-based ***IndexSieveOfAtkin***, which juggles between numbers and indices. But on the other hand ***SieveOfAtkin*** performs Euclidian divisions by 12, whereas ***IndexSieveOfAtkin*** does divisions by 6. This is due to the conversion of number $n$ into its index $k = (n - 3)/2$. Furthermore, the latter only performs the sieve on odd numbers, which means effectively the memory space for the sieve is twice smaller.

On the graph 3-3 below, we plot the computation time in seconds for both algorithms. The curve $T_3$ corresponds to ***SieveOfAtkin*** and the curve $T_4$ to ***IndexSieveOfAtkin***. We observe empirically that computation time of both algorithms looks slightly higher than linear, even though theoretically the number of operations appears to be linear in $N_{Max}$. Details of the Maple options used to get the regression are given in appendix 8.2.

**Graph 3-3**: computation time $T(N_{Max})$ in seconds for both algorithms (Sieve of Atkin)

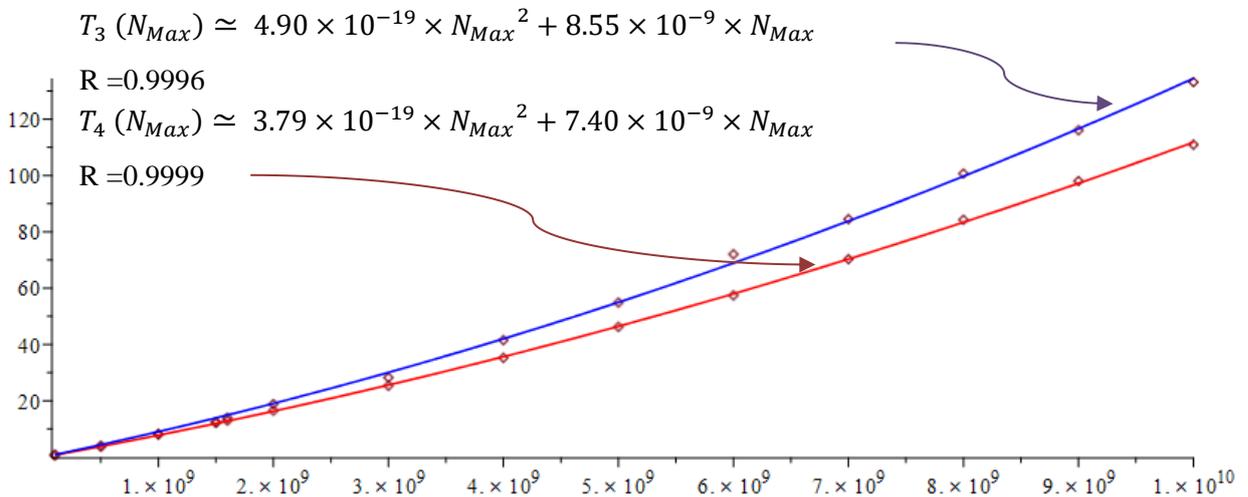

$$T_3(N_{Max}) \simeq 4.90 \times 10^{-19} \times N_{Max}^2 + 8.55 \times 10^{-9} \times N_{Max}$$
R =0.9996
$$T_4(N_{Max}) \simeq 3.79 \times 10^{-19} \times N_{Max}^2 + 7.40 \times 10^{-9} \times N_{Max}$$
R =0.9999

The second algorithm is faster for larger values of $N_{Max}$, roughly for $N_{Max} > 10^9$. For such values the cost of encoding numbers to indices is offset by the gain on modulo operations and halving the size of the sieve. We note also that memory size is halved for the second algorithm.

# 4 Wheel sieve with indices

We first describe Pritchard's wheel sieve. Then we adapt it to indices and discuss a way to generate the integers of the *turning wheel*.

## 4.1 Description of Pritchard's wheel sieve

This description is based on [7] and [4]. The wheel sieve operates by generating a set of numbers that are coprime with the first $q$ prime numbers. The second of these is the next prime, multiples of which are then eliminated (*by turning the wheel*).

More precisely, let $p_0 = 2, p_1 = 3 \dots$ the sequence of prime numbers and let:

$$\Pi_q = \prod_{k=0}^{q} p_k$$

$$\mathcal{R}(m) = \{x \in [\![1, m-1]\!] \mid \gcd(x, m) = 1\}$$

$$\mathcal{W}_q = \mathcal{R}(\Pi_q)$$

The following proposition describes a "turn of the wheel".

**<u>Proposition 4-1-1:</u>** We have the following inductive formula for $\mathcal{W}_q$:

$$\mathcal{W}_0 = \{1\}, \mathcal{W}_1 = \{1,5\}, \mathcal{W}_2 = \{1,7,11,13,17,19,23,29\}$$

$$\forall q \in \mathbb{N}, \mathcal{W}_{q+1} = \left( \bigcup_{x=0}^{p_{q+1}-1} (\mathcal{W}_q + x\Pi_q) \right) \setminus p_{q+1}[\![1, \Pi_q - 1]\!]$$

<u>*Proof:*</u> The Chinese theorem ensures that $m \in \mathcal{W}_{q+1}$ if and only if $m \bmod \Pi_q \in \mathcal{W}_q$ and $m \notin p_{q+1}\mathbb{N}$. This gives the desired set equality.

Furthermore, induction formula for $\mathcal{W}_q$ can also be used to recursively build the sequence of prime numbers:

**<u>Proposition 4-1-2:</u>** The second smallest element of $\mathcal{W}_q$ ($q \geq 1$) is the next prime $p_{q+1}$.

<u>*Proof:*</u> The first element is 1, which is obviously not prime. For $q \geq 1$, $p_q \geq 3$ and from proposition 4-1-1 we can show (see corollary 4-2-2 later on) that $\mathcal{W}_q$ has at least 2 elements. The second one must then be the smallest integer coprime with $p_0 \dots p_q$, and thus must be $p_{q+1}$.

The elements of $\mathcal{W}_q$ are called pseudo-primes (at order $q$). Some of them are primes and others are not. However, we have a boundary condition to identify some of the primes:

**<u>Proposition 4-1-3:</u>** All integers in $\mathcal{W}_q$ and less than $p_q^2$ are sure to be primes.

<u>*Proof:*</u> Any integer less than $p_q^2$ is either prime or has a divisor among $p_0 \dots p_q$. The latter is impossible by definition of $\mathcal{W}_q$.

To enumerate primes up to $N_{Max}$, we thus have to keep turning the wheel as long as $p_{q+1}^2 < N_{Max}$.

As $\Pi_q$ grows exponentially (in particular it can be easily proven from Bertrand's postulate that $\Pi_q > p_q^2$ from $q = 2$), while we are only interested in pseudo-primes up to $N_{Max}$, we may replace in practice $\mathcal{W}_q$ by $\mathcal{W}_q^{N_{Max}} = \mathcal{W}_q \cap [\![1, N_{Max}]\!]$.

**Proposition 4-1-4:** The following inductive formula (or wheel turn) is true for all $N_{Max}$:

$$\forall q \in \mathbb{N}, \mathcal{W}_{q+1}^{N_{Max}} = \left[\left(\bigcup_{x=0}^{\max\left(p_{q+1}-1, \left\lfloor \frac{N_{Max}}{\Pi_q}\right\rfloor\right)} (\mathcal{W}_q + x\Pi_q)\right) \setminus p_{q+1}\left[\!\left[1, \left\lfloor \frac{N_{Max}}{p_{q+1}}\right\rfloor\right]\!\right]\right] \cap [\![1, N_{Max}]\!].$$

Furthermore, if $N_{Max} > 9$, then as soon as $p_q^2 \geq N_{Max}$, $P_{N_{Max}} = \{p_0 \ldots p_q\} \cup \left(W_q^{N_{Max}} \setminus \{1\}\right)$.

*Proof:* By double inclusion (cf. proof of proposition 4-2-3). The second identity comes from the fact that if $N_{Max} > 9$, $p_q^2 \geq N_{Max}$ implies $q \geq 2$.

Thus, when we turn the wheel, we remove integers that are, for a given $m \in \mathcal{W}_q$, and $x, y$ integers, of the form:

$$m + x\Pi_q = yp_{q+1}$$

One way to do that is to remove all multiples of $p_{q+1}$. We will show however in section 4.2 that there is a relationship between the value of $x$, the multiples of $\Pi_q$ which are added to $\mathcal{W}_q$, and the composite numbers $yp_{q+1}$ which must be removed of the wheel $\mathcal{W}_{q+1}$, so that the index $x$ to remove can be predicted from $m$ or conversely.

### 4.2 Index wheel sieve

**Definition 4-2**: We note $\Pi'_q$ the product of all odd primes up to $p_q$, i.e. $\Pi_q = 2\Pi'_q$.

We also note:

$$N(m, a, q) = m\Pi_q + a$$

and, with $a'$ the index of $a$:

$$k(m, a', q) = \frac{N(m, 2a' + 3, q) - 3}{2} = m\Pi'_q + a'$$

the index of $N(m, a, q)$.

We let $\mathcal{W}'_q$ be the set of indices corresponding to $\mathcal{W}_q$, with 1 replaced by $\Pi_q + 1$ (which index is $\Pi'_q - 1$):

$$\mathcal{W}'_q = \left\{\frac{n-3}{2}, n \in \mathcal{W}_q \setminus \{1\}\right\} \cup \{\Pi'_q - 1\}$$

In this section, we describe how we adapt the wheel sieve to work with indices of odd integers. The limit $N_{Max}$ is supposed to be an odd integer of index $k_{Max}$.

*Recurrence relation verified by the index wheel sieve:*

The initial index wheels are $\mathcal{W}'_0 = \{0\}$, $\mathcal{W}'_1 = \{1,2\}$, $\mathcal{W}'_2 = \{2,4,5,7,8,10,13,14\}$.

**Remark 4-2-1:** The first element of $\mathcal{W}'_q$ is the index of the prime number $p_{q+1}$. $\mathcal{W}'_q$ is included in $\left[\!\left[\frac{p_{q+1}-3}{2}, \Pi'_q - 1\right]\!\right]$.

*Proof*: Since we remapped 1 to $\Pi_q + 1$ in $\mathcal{W}_q$ to define $\mathcal{W}'_q$, and because the indexing map is increasing, the first element of $\mathcal{W}'_q$ is the index of prime $p_{q+1}$ from proposition 4-1-2 (we note that it works even for $q = 0$), and its last element is $\Pi'_q - 1$.

**Proposition 4-2-1:** The <u>index wheel sieve</u> is the only sequence of sets verifying:

$$\mathcal{W}'_0 = \{0\}$$

$$\forall q \in \mathbb{N}, \mathcal{W}'_{q+1} = \left( \bigcup_{m=0}^{p_{q+1}-1} \left(\mathcal{W}'_q + m\Pi'_q\right) \right) \setminus \left\{\frac{p_{q+1}-3}{2} + y'p_{q+1}, y' \in \left[\!\left[0, \Pi'_q - 1\right]\!\right]\right\}$$

Furthermore, indices in the wheel $\mathcal{W}'_q$ up to $k_{Max}$ correspond to all remaining prime numbers up to $N_{Max}$ (on top of $p_0 \ldots p_q$) as soon as:

$$\frac{{p_q}^2 - 3}{2} \geq k_{Max}$$

*Proof*: This comes from the definition 4-2 of the index wheel sieve, the proposition 4-1-1 and from observing that the index of any odd multiple $yp_q$ of $p_q$ is of the form:

$$\frac{yp_q - 3}{2} = \frac{p_q - 3}{2} + y'p_q, y' = \frac{y-1}{2}$$

If we let $p = 2i + 3$, this corresponds to the definition of $k(y', i)$ in [8]: $k(y', i) = i + (2i + 3)y'$.

*Eliminating multiples of the next prime by solving a Diophantine equation*:

**Proposition 4-2-2:** For a given $c \in \left[\!\left[0, \Pi'_q - 1\right]\!\right]$, there exists a unique $(m_c, y_c) \in \left[\!\left[0, p_{q+1}-1\right]\!\right] \times \mathcal{W}_q$ such that $c + m_c\Pi'_q = y_c p_{q+1}$. Furthermore, $m_c$ only depends of $c \bmod p_{q+1}$, $m_0 = 0$ and for $c_1 = \left(-\Pi'_q\right) \bmod p_{q+1}$,

$$m_{c_1} = 1.$$

For all $c \in \mathcal{W}_q$ one has $c \bmod p_{q+1} = m_c c_1 \bmod p_{q+1}$

**Remark 4-2-2:** Using indices, we must solve $(m, y')$ in the following equations for $a' \in \mathcal{W}'_q$:

$$a' + m\Pi'_q = \frac{p_{q+1}-3}{2} + y'p_{q+1}$$

so we will let $c = a' - \frac{p_{q+1}-3}{2}$.

*Proof*: Because $\Pi'_q$ and $p_{q+1}$ are coprime, existence and unicity of the solution are well-known. In [9] we introduced the concept of normalizer of such a Diophantine equation, and have shown its additive and multiplicative property.

Clearly if $c \equiv d \ [p_{q+1}]$ then $(m_c - m_d)\Pi'_q \equiv 0 \ [p_{q+1}]$ and as $\Pi'_q$ and $p_{q+1}$ are coprime, $m_c \equiv m_d \ [p_{q+1}]$.

Also, because $0 + 0.\Pi'_q = 0.p_{q+1}$ we deduce that $m_0 = 0$.

Then from the fact that $c_1 + \Pi'_q \in p_{q+1}\mathbb{Z}$ we get that $m_{c_1} = 1$.

Furthermore, for all $c$, by multiplicative property:

$$m_{m_c c_1} \equiv m_c . m_{c_1} \equiv m_c \ [p_{q+1}]$$

Thus, $c \equiv -m_c \Pi'_q \equiv -m_{m_c c_1}\Pi'_q \equiv m_c c_1 \ [p_{q+1}]$.

This proposition gives us an effective way of building all couples $(c, m_c)$ modulo $p_{q+1}$: start from $(c_1, 1)$ and add it to itself (modulo $p_{q+1}$) up to $p_{q+1} - 1$ times (the last time we will get the couple $(0, 0 = m_0)$).

**Corollary 4-2-2:** $\mathcal{W}_q$ and $\mathcal{W}'_q$ have $\prod_{k=1}^{q}(p_k - 1)$ elements.

*Proof*: Let us proceed by induction on $q$. The property is true for $q = 0$. Assume it is true for a given $q \in \mathbb{N}$. From proposition 4-2-1,

$$\mathcal{W}'_{q+1} = \left( \bigcup_{m=0}^{p_{q+1}-1} \left( \mathcal{W}'_q + m\Pi'_q \right) \right) \setminus \left\{ \frac{p_{q+1}-3}{2} + y'p_{q+1}, y' \in [\![0, \Pi'_q - 1]\!] \right\}.$$

Thus $\bigcup_{m=0}^{p_{q+1}-1}(\mathcal{W}'_q + m\Pi'_q) = \bigcup_{c' \in \mathcal{W}'_q}(c' + \Pi'_q[\![0, p_{q+1} - 1]\!])$ has exactly $p_{q+1}\prod_{k=1}^{q}(p_k - 1)$ elements, from which we must remove the indices of multiples of $p_{q+1}$. For a given $c' \in \mathcal{W}'_q$, from proposition 4-2-2 there is exactly one couple $(m, y)$ such that:

$$c' + m\Pi'_q = \frac{p_{q+1}-3}{2} + y'p_{q+1}$$

i.e. there is only one element of $c' + \Pi'_q[\![0, p_{q+1} - 1]\!]$ in $\left\{ \frac{p_{q+1}-3}{2} + y'p_{q+1}, y' \in [\![0, \Pi'_q - 1]\!] \right\}$. So in total there are exactly $\prod_{k=1}^{q}(p_k - 1)$ elements in $\left( \bigcup_{m=0}^{p_{q+1}-1}(\mathcal{W}'_q + m\Pi'_q) \right) \cap \left\{ \frac{p_{q+1}-3}{2} + y'p_{q+1}, y' \in [\![0, \Pi'_q - 1]\!] \right\}$, thus $(p_{q+1} - 1)\prod_{k=1}^{q}(p_k - 1) = \prod_{k=1}^{q+1}(p_k - 1)$ elements in $\mathcal{W}'_{q+1}$.

**Proposition 4-2-3:** $\mathcal{W}'_q{}^{k_{Max}} = \mathcal{W}'_q \cap [\![0, k_{Max}]\!]$ verifies the following induction property.

For all $q \in \mathbb{N}$, $\mathcal{W}'_{q+1}{}^{k_{Max}}$ is equal to:

$$\left( \left( \bigcup_{m=0}^{\min\left(p_{q+1}-1, \left\lfloor \frac{k_{Max}}{\Pi'_q} \right\rfloor\right)} \left( \mathcal{W}'_q{}^{k_{Max}} + m\Pi'_q \right) \right) \setminus \left\{ \frac{p_{q+1}-3}{2} + y'p_{q+1}, y' \in \left[\!\!\left[ 0, \min\left( \Pi'_q - 1, \left\lfloor \frac{2k_{Max}+3}{2p_{q+1}} - \frac{1}{2} \right\rfloor \right) \right]\!\!\right] \right\} \right) \cap [\![0, k_{Max}]\!]$$

<u>*Proof*</u>: Let $x \in \mathcal{W}'_{q+1}{}^{k_{Max}}$. From proposition 4-2-1, there exists $c' \in \mathcal{W}'_q$, $m \in [\![0, p_{q+1} - 1]\!]$ such that $x = c' + m\Pi'_q$. But $x \leq k_{Max}$ so $m \leq \lfloor k_{Max}/\Pi'_q \rfloor$. Furthermore, $x \notin \left\{ \frac{p_{q+1}-3}{2} + y'p_{q+1}, y' \in [\![0, \Pi'_q - 1]\!] \right\}$ so a fortiori:

$$x \notin \left\{ \frac{p_{q+1}-3}{2} + y'p_{q+1}, y' \in \left[\!\!\left[ 0, \min\left( \Pi'_q - 1, \left\lfloor \frac{2k_{Max}+3}{2p_{q+1}} - \frac{1}{2} \right\rfloor \right) \right]\!\!\right] \right\}.$$

Conversely, let $x \in \left( \bigcup_{m=0}^{\min(p_{q+1}-1, \lfloor k_{Max}/\Pi'_q \rfloor)} \left( \mathcal{W}'_q{}^{k_{Max}} + m\Pi'_q \right) \right) \cap [\![0, k_{Max}]\!]$ such that $x \notin \left\{ \frac{p_{q+1}-3}{2} + y'p_{q+1}, y' \in \left[\!\!\left[ 0, \min\left( \Pi'_q - 1, \left\lfloor \frac{2k_{Max}-3}{2p_{q+1}} - \frac{1}{2} \right\rfloor \right) \right]\!\!\right] \right\}$. The first condition means that $x \in \mathcal{W}'_{q+1}$ if $x \notin \left\{ \frac{p_{q+1}-3}{2} + y'p_{q+1}, y \in [\![0, \Pi'_q - 1]\!] \right\}$. But if that were the case, there would be $y' \in [\![1, \Pi'_q - 1]\!]$ such that $x = \frac{p_{q+1}-3}{2} + y'p_{q+1}$. Thus $y \leq \frac{k_{Max}-(p_{q+1}-3)/2}{p_{q+1}} = \frac{2k_{Max}+3}{2p_{q+1}} - \frac{1}{2}$, which cannot happen because $x \notin \left\{ \frac{p_{q+1}-3}{2} + y'p_{q+1}, y' \in \left[\!\!\left[ 0, \min\left( \Pi'_q - 1, \left\lfloor \frac{2k_{Max}+3}{2p_{q+1}} - \frac{1}{2} \right\rfloor \right) \right]\!\!\right] \right\}$.

### 4.3 Wheel sieve algorithms

As per sections 4.1 and 4.2, the wheel sieve algorithms will consist in two steps:

(A) A first step where the wheel will always grow, as long as $\Pi_q < N_{Max}$, or:
$$\Pi'_q - 1 < k_{Max},$$

(B) A second step where we will no longer grow the wheel, but will have to keep eliminating composite numbers, as long as $p_q^2 < N_{Max}$, or:
$$\frac{p_q^2 - 3}{2} < k_{Max}.$$

This is equivalent to saying that we replace $\mathcal{W}_{q+1}$ by $\mathcal{W}_{q+1}^{N_{Max}}$ and similarly $\mathcal{W}'_{q+1}$ by $\mathcal{W}'_{q+1}{}^{k_{Max}}$. During step (B) we do not add new pseudo-primes, only remove those that we rule out as multiples of the next prime. Because $\Pi_q$ grows exponentially, there will generally be more iterations in step (B) than in step (A).

<u>Quick description of the steps of the index wheel sieve algorithm (see appendix for full algorithm):</u>

As for the previous algorithms, we note $L_p$ the list of primes and $i_l$ its number of elements. $IL_p$ represents the list of indices of odd primes, and $SIL_p$ the list of indices of squared odd primes. At step $q$, $L_p$ will contain all primes up to $p_q^2$, coming from the wheel $\mathcal{W}'_q$, $IL_p$ and $SIL_p$ being filled with the corresponding indices.

1- Intialisation of the sieve for $q = 1$: $L_p = \{2,3,5,7\}$, $i_l = 4$ $IL_p = \{0,1,2\}$, $SIL_p = \{3,11,23\}$ and $\mathcal{W}'_1 = \{1, 2\}$ with $\Pi'_1 = 3$.
2- While $\Pi'_q < k_{max}$ (step A):
    a. We take $p_{q+1}$ from $L_p$ (or equivalently the first element of $\mathcal{W}_q$). The list of pairs $(c, m_c)$ such that $c + m_c\Pi'_q$ has to be eliminated is then computed, according to proposition 4-2-2. Then we build the wheel $\mathcal{W}'_{q+1}$.
    b. Once this is done primes in the interval $[\![p_{i_{l-1}} + 2, p_{q+1}^2 - 2]\!]$ are added to $L_p$ and $i_l$, $IL_p$ and $SIL_p$ are updated accordingly. Indices of the primes to add are those in

$$\mathcal{W}'_{q+1} \cap [\![IL_p(i_l - 2) + 1, SIL_p(q) - 1]\!].$$

3- While $SIL_p(q) < k_{max}$ (step B):
   a. Remove indices of multiples of $p_{q+1}$ from $\mathcal{W}'_q{}^{k_{Max}}$ to get $\mathcal{W}'_{q+1}{}^{k_{Max}}$.
   b. Once this is done primes in the interval $[\![p_{i_{l-1}} + 2, p_{q+1}^2 - 2]\!]$ are added to $L_p$ and $i_l$, $IL_p$ and $SIL_p$ are updated accordingly. Indices of the primes to add are those in $\mathcal{W}'_{q+1} \cap [\![IL_p(i_l - 2) + 1, SIL_p(q) - 1]\!]$.

**Remark 4-3-1:** Let $k_1$ and $k_2$ be the indices of two odd numbers, respectively $n_1$ and $n_2$, such as $n_2 - n_1 > 0$. Let $\alpha = k_2 - k_1$. The difference between the indices $n_1^2$ and $n_2^2$ is:

$$\beta = 2\alpha^2 + 2\alpha n_1.$$

Furthermore, if $m$ is another integer, the difference between the indices of $n_1 m$ and $n_2 m$ is:

$$\gamma = \alpha m.$$

*Proof*: Note that $n_2 - n_1 = 2\alpha$ and thus:

$$\frac{n_2^2 - 3}{2} - \frac{n_1^2 - 3}{2} = \frac{1}{2}(n_2 - n_1)(n_2 + n_1) = \alpha(n_2 + n_1) = \alpha(2n_1 + 2\alpha) = \beta.$$

Similarly:

$$\frac{n_2 m - 3}{2} - \frac{n_1 m - 3}{2} = \alpha m = \gamma.$$

This last remark is used in steps 2-b. and 3-b. to fill $SIL_p$ and to perform step 3-a.

**Remark 4-3-2:** The index wheel sieve involves operations with reduced input size compared with the number version. This is clear from remark 4-3-1 where $\beta$ is exactly half of $n_2^2 - n_1^2$, for instance. Similarly $\Pi'_q$ is half of $\Pi_q$ so modulo operation input is also reduced.

## 4.4 Performance of algorithms

In this section, we present results from the previous algorithm of index wheel sieve, which we compare with a similar one on numbers (unspecified for to avoid a lengthy duplication). These results are similar to those obtained in the previous sections. As for the sieve of Atkin, we did not go for refinements that give a better time complexity, so theoretical complexity in terms of number of operations is $O(N)$ for both algorithms.

On the graph 4-4 below, we plot the computation time in seconds for both algorithms, for $N_{Max}$ up to $6.10^9$. The curve $T_5$ corresponds to the the algorithm ***WheelSieveReference*** and the curve $T_6$ corresponds to the the algorithm ***IndexWheelSieve***. The correlation coefficient $R$ of each regression is given on the graph. Details of the Maple options used to get the regression are given in appendix 8.3. We notice that complexity of both algorithms again seems empirically slightly higher than linear.

**Graph 4-4**: computation time $T(N_{Max})$ in seconds for both algorithms (Wheel sieve)

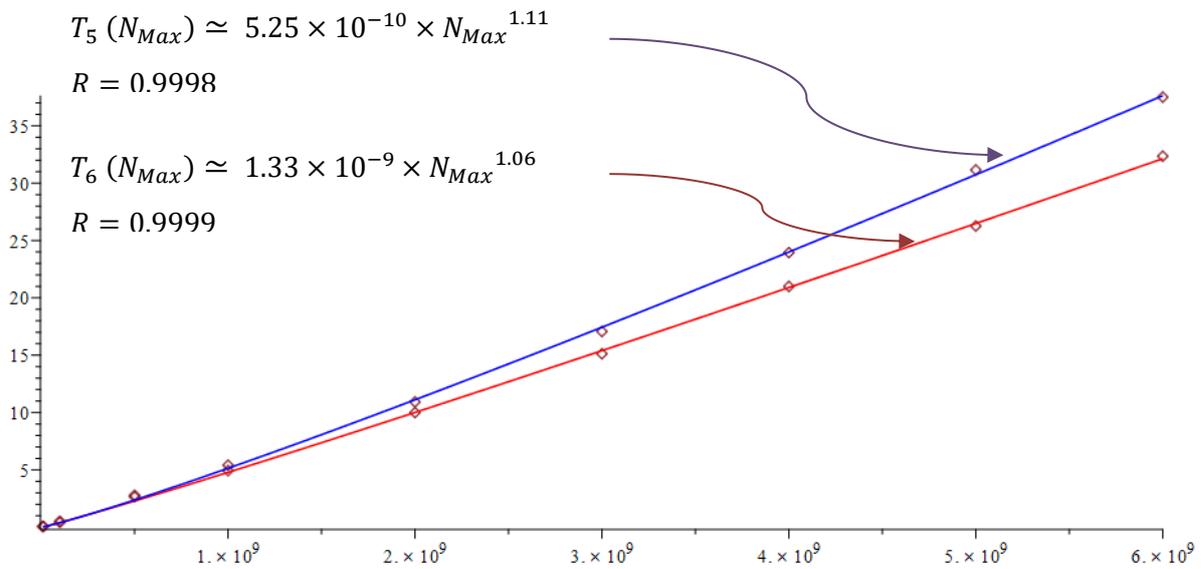

$T_5(N_{Max}) \simeq 5.25 \times 10^{-10} \times N_{Max}^{1.11}$

$R = 0.9998$

$T_6(N_{Max}) \simeq 1.33 \times 10^{-9} \times N_{Max}^{1.06}$

$R = 0.9999$

Complexity is reduced by using indices, due to reduction of input size in the modulo and the multiplication operations (see Remark 4-3-2) and despite a higher number of operations with the algorithm **IndexWheelSieve**. Moreover, the amount of memory space used with indices is halved, due to the fact that we avoid even numbers completely.

# 5 Conclusion

In theory, indices are a way to work with odd numbers only by not representing even numbers. Most mathematical relations must be reformulated for indices, which lead to a higher number of (conversion) operations, but in return the input size of other operations is reduced. In this article, we have shown how this indexing translates into optimized algorithms in applied mathematics. From a basic primality test implementation, to the sieve of Atkin and Pritchard's wheel sieve, indices speeded up these algorithms, not by changing their complexity but by reducing the time cost by a constant factor, and generally also made them more efficient from a memory point of view.

Acknowledgments: We would like to thank François-Xavier VILLEMIN for his attentive comments and suggestions.

## 7 APPENDIX: ALGORITHM OF THE INDEX WHEEL SIEVE

This algorithm enumerates odd primes up to the limit $N_{Max}$. It is composed of a main function that is called **IndexWheelSieve** and the following auxilliary other functions:

7-2- **DiophantineSolutions**$(Prime, \Pi'_q)$

7-3- **WheelTurn**$(\mathcal{W}'_q, q, Prime, PrimeIndex, \Pi'_q, k_{Max})$

7-4- **RemoveMultiples**$(SquarePrimeIndex, Prime, \mathcal{W}'_q)$

7-5- **GetNewPrimes**$(\mathcal{W}'_q, q, L_p, i_l, IL_p, SIL_p)$

Some marginal optimizations can still be performed, for instance modulo operations inside a loop can be replaced by substractions, and memory can be managed better. For the sake of readability we leave these optimizations out of scope.

---

**Algorithm 7-1 *IndexWheelSieve*$(N_{Max})$:** $N_{Max} > 9$ is an odd integer.

This function returns the list of all prime numbers up to $N_{Max}$.

---

***First step : intialisation of variables***

$L_p \leftarrow \{2,3,5,7\}$        → Dynamic list of primes

$i_l \leftarrow 4$        → Number of primes in the list

$k_{Max} \leftarrow (N_{Max} - 3)/2$        → Index of $N_{Max}$

$IL_p \leftarrow \{0,1,2\}$

$SIL_p \leftarrow \{3, 11, 23\}$

$\mathcal{W}'_q \leftarrow \{1, 2\}$

$\Pi'_q \leftarrow 3$

$q \leftarrow 1$

***Second step : Wheel inflation***.

**Do**

  $Prime \leftarrow L_p(q + 1)$

  $PrimeIndex \leftarrow IL_p(q)$

  $\Pi'_{q+1} \leftarrow \Pi'_q \times Prime$

  $\rightarrow$ Compute values of the new wheel from the previous one

  $\mathcal{W}'_q \leftarrow$ **WheelTurn**$(\mathcal{W}'_q, q, Prime, PrimeIndex, \Pi'_q, k_{Max})$

  **GetNewPrimes**$(\mathcal{W}'_q, q, L_p, i_l, IL_p, SIL_p)$

  $\Pi'_q \leftarrow \Pi'_{q+1}$

  $q \leftarrow q + 1$

**While** $k_{Max} > \Pi'_q$

***Third step : Wheel deflation.***

**While**   $SIL_p(q - 1) < k_{Max}$

  $Prime \leftarrow L_p(q + 1)$

  $SquarePrimeIndex \leftarrow SIL_p(q)$

  $\mathcal{W}'_q \leftarrow$ **RemoveMultiples**$(SquarePrimeIndex, Prime, \mathcal{W}'_q)$

  **GetNewPrimes** $(\mathcal{W}'_q, q, L_p, i_l, IL_p, SIL_p)$

  $q \leftarrow q + 1$

**End While**

**Return** $(L_p, i_l)$

---

---

**Algorithm 7-2 DiophantineSolutions**$(Prime, \Pi'_q)$

---

$c_1 \leftarrow Prime - \left( \Pi'_q \bmod Prime \right)$          $\rightarrow$ Solution such that $m = 1$

$c \leftarrow 0$

$Solutions \leftarrow \{0 \dots 0\}$       $\rightarrow$ Array of size $Prime$

**For** $m = 1$ **To** $Prime - 1$ **Do**

 $c \leftarrow (c + c_1) \bmod Prime$

 $Solutions(c) \leftarrow m$

**End For**

**Return** $Solutions$

---

---

**Algorithm 7-3 WheelTurn**$(\mathcal{W}'_q, q, Prime, PrimeIndex, \Pi'_q, k_{Max})$

This function computes $\mathcal{W}'_{q+1}$ by duplicating the wheel $\mathcal{W}'_q$ and removing indices of multiples of $Prime = p_{q+1}$.

---

***First step : Compute all the pairs*** $(c, m_c)$ *in the function* **DiophantineSolutions**

$Solutions \leftarrow$ **DiophantineSolutions**$(Prime, \Pi'_q)$

***Second step : Iteration***

$WheelSize \leftarrow$ **Size**$(\mathcal{W}'_q)$

$Table \leftarrow$ **Range**$(\{\}, Prime)$

**For** $j = 0$ **To** $WheelSize - 1$ **Do**

 $a' \leftarrow \mathcal{W}'_q(j)$

 $c \leftarrow (a' - PrimeIndex) \bmod Prime$

 $m \leftarrow Solutions(c)$

 **For** $y = 0$ **To** $PrimeNumber - 1$ **Do**

  $n \leftarrow a' + y\Pi'_q$

  **If** $n > k_{max}$ **Do**

   **Break**

  **End If**

  **If** $y \neq m$ **Do**

   **Append**$(Table(y), n)$

  **End If**

 **End For**

**End For**

***Third step : Build $\mathcal{W}'_{q+1}$ by concatenation***

$\mathcal{W}'_{q+1} \leftarrow \{\}$

**For** $y = 0$ **To** $PrimeNumber - 1$ **Do**

  **Concatenate**$(\mathcal{W}'_{q+1}, Table(y))$

**End For**

**Return** $\mathcal{W}'_{q+1}$

---

**Algorithm 7-4 RemoveMultiples**$(SquarePrimeIndex, Prime, \mathcal{W}'_q)$

---

$\mathcal{W}'_{q+1} \leftarrow \{\}$

$NextMultiple \leftarrow SquarePrimeIndex$

**For** $j = 1$ **To Size**$(\mathcal{W}'_q) - 1$ **Do**

  **If** $\mathcal{W}'_q(j) > NextMultiple$ **Do**

   $NextMultiple \leftarrow NextMultiple + Prime$

   $j \leftarrow j - 1$

  **Else If** $\mathcal{W}'_q(j) = NextMultiple$ **Do**

   $NextMultiple \leftarrow NextMultiple + Prime$

  **Else**

    **Append**$(\mathcal{W}'_{q+1}, \mathcal{W}'_q(j))$

  **End If**

**End For**

**Return** $\mathcal{W}'_{q+1}$

---

**Algorithm 7-5 GetNewPrimes**$(\mathcal{W}'_q, q, L_p, i_l, IL_p, SIL_p)$

This function adds new primes to the list $L_p$ and updates $i_l$ and the other lists $IL_p$ and $SIL_p$ (all passed by reference).

---

$SquareIndex \leftarrow SIL_p(q + 1)$

$j \leftarrow i_l - q - 2$                 → Offset to take into account already known primes

$NewPrimeIndex \leftarrow \mathcal{W}'_q(j)$

**While** $NewPrimeIndex < SquareIndex$ **Do**

   $IL_p(i_l - 1) \leftarrow NewPrimeIndex$

   $\alpha \leftarrow IL_p(i_l - 1) - IL_p(i_l - 2)$

   $SIL_p(i_l - 1) \leftarrow SIL_p(i_l - 2) + 2\alpha^2 + 2\alpha L_p(i_l - 1)$

   $L_p(i_l) \leftarrow L_p(i_l - 1) + 2\alpha$

   $i_l \leftarrow i_l + 1$

   $j \leftarrow j + 1$

   $NewPrimeIndex \leftarrow \mathcal{W}'_q(j)$

**End While**

---

<h1 style="text-align:center;">8 APPENDIXES: MAPLE REGRESSIONS</h1>

Here are numeric values obtained from our implementation (Visual Studio C++ 2012) of the algorithms presented in this article.

<h2 style="text-align:center;">8.1 BASIC PRIMALITY TEST AND PRIMES ENUMERATION</h2>

In table 8.1, numeric values of $T_1(N_{Max})$ and $T_2(N_{Max})$ are obtained respectively from the **PrimeEnumeration** and **IndexPrimeEnumeration** algorithms.

<p style="text-align:center;"><strong>Table 8.1:</strong> numeric values of $T_1(N_{Max})$ and $T_2(N_{Max})$ in seconds.</p>

| $N_{Max}$ | $10^7$ | $10^8$ | $5 \times 10^8$ | $10^9$ | $2 \times 10^9$ | $3 \times 10^9$ | $4 \times 10^9$ |
|---|---|---|---|---|---|---|---|
| $T_1(N_{Max})$ | 2.403 | 56.031 | 493.163 | 1306.884 | 3414.713 | 6271.249 | 8908.814 |
| $T_2(N_{Max})$ | 2.375 | 54.725 | 487.568 | 1275.921 | 3329.573 | 6105.386 | 8664.438 |

To fit these observations, Maple's **NonlinearFit** function is used with the parameters below. Initial values for $a$ and $b$ were determined empirically.

```
NonlinearFit(a × n^b / ln(n), X, Y, n, initialvalues = [a = 5.9 × 10⁻⁹, b = 1.41],
             output = [leastsquaresfunction, residuals])
```

We get the following mathematical relationships:

$$T_1(N_{Max}) \simeq 5.79409775129480 \times 10^{-9} \times \frac{n^{1.40966993452829}}{\ln(n)}, \; R = .99962000$$

$$T_2(N_{Max}) \simeq 6.10602965467609 \times 10^{-9} \times \frac{n^{1.406040046365699}}{ln(n)}, \; R = .99962009$$

## 8.2 THE SIEVE OF ATKIN

In table 8.2, numeric values of $T_3(N_{Max})$ and $T_4(N_{Max})$ are obtained respectively from the **SieveOfAtkin** and **IndexSieveOfAtkin** algorithms.

| $N_{Max}$ | $10^8$ | $5 \times 10^8$ | $10^9$ | $1.5 \times 10^9$ | $1.6 \times 10^9$ | $2 \times 10^9$ | $3 \times 10^9$ |
|---|---|---|---|---|---|---|---|
| $T_3(N_{Max})$ | 0.719 | 3.797 | 8.033 | 12.48 | 13.967 | 18.843 | 28.217 |
| $T_4(N_{Max})$ | 0.727 | 3.921 | 8.225 | 12.152 | 12.953 | 16.507 | 25.342 |

| $N_{Max}$ | $4 \times 10^9$ | $5 \times 10^9$ | $6 \times 10^9$ | $7 \times 10^9$ | $8 \times 10^9$ | $9 \times 10^9$ | $10^{10}$ |
|---|---|---|---|---|---|---|---|
| $T_3(N_{Max})$ | 41.534 | 54.871 | 72.044 | 84.511 | 100.727 | 116.093 | 133.184 |
| $T_4(N_{Max})$ | 35.27 | 46.261 | 57.418 | 70.311 | 84.291 | 98.047 | 110.96 |

This time we used Maple's function **Fit** as below:

$$\texttt{Fit}(a \times n^2 + b \times n, \ \texttt{X, Y}, \ n, \ \texttt{summarize = embed})$$

We get the following mathematical relationships:

$T_3(N_{Max}) \simeq 4.90268369826396 \times 10^{-19} \times {N_{Max}}^2 + 8.54576412559177 \times 10^{-9} \times N_{Max}, \quad R = .999647$

$T_4(N_{Max}) \simeq 3.78795281632082 \times 10^{-19} \times {N_{Max}}^2 + 7.39595089422000 \times 10^{-9} \times N_{Max},$

$R = .999926$

## 8.3 WHEEL SIEVE WITH INDICES

In table 8.3, numeric values of $T_5(N_{Max})$ and $T_6(N_{Max})$ are obtained respectively from the **WheelSieveReference** and **IndexWheelSieve** algorithms.

| $N_{Max}$ | $10^7$ | $10^8$ | $5 \times 10^8$ | $10^9$ | $2 \times 10^9$ | $3 \times 10^9$ | $4 \times 10^9$ | $5 \times 10^9$ | $6 \times 10^9$ |
|---|---|---|---|---|---|---|---|---|---|
| $T_5(N_{Max})$ | 0.071 | 0.496 | 2.783 | 5.407 | 10.931 | 17.070 | 23.944 | 31.150 | 37.501 |
| $T_6(N_{Max})$ | 0.064 | 0.457 | 2.657 | 4.936 | 9.995 | 15.121 | 20.995 | 26.260 | 32.351 |

We used again **NonlinearFit** with empirically determined initial values $a$ and $b$:

$\texttt{NonlinearFit}(a \times n^b, \ \texttt{X, Y}, \ n, \ \texttt{initialvalues} = [a = 1.97461115539853 \times 10^{-6},$
$\quad\quad b = 1.1], \ \texttt{output} = [\texttt{leastsquaresfunction, residuals}]).$

We get the following mathematical relationships:

$T_5(N_{Max}) \simeq 5.25118782575365 \times 10^{-10} \times n^{1.11016647384427}, R = .99982444$

$T_6(N_{Max}) \simeq 1.33020583039257 \times 10^{-9} \times n^{1.06187203820827}, \ R = .99986693$